\documentclass[
11pt,%
tightenlines,%
twoside,%
onecolumn,%
nofloats,%
nobibnotes,%
nofootinbib,%
superscriptaddress,%
noshowpacs,%
centertags]%
{revtex4-2}

\usepackage{amssymb}
\usepackage{amsfonts}
\usepackage[tbtags]{amsmath}
\usepackage{amscd}
\usepackage{amsthm}
\usepackage{mathtext}
\usepackage[T2A]{fontenc}

\usepackage[cp1251]{inputenc}
\usepackage[russian,english]{babel}
\usepackage[all]{xy}

\newcommand{\Raise}{\oper{P}}
\newcommand{\Lower}{\oper{Q}}

\usepackage{amsmath}

\usepackage{mathtools}

\makeatletter
\newcommand{\hathat}[1]{%
\begingroup%
  \let\macc@kerna\z@%
  \let\macc@kernb\z@%
  \let\macc@nucleus\@empty%
  \hat{\mathchoice%
    {\raisebox{.2ex}{\vphantom{\ensuremath{\displaystyle #1}}}}%
    {\raisebox{.2ex}{\vphantom{\ensuremath{\textstyle #1}}}}%
    {\raisebox{.16ex}{\vphantom{\ensuremath{\scriptstyle #1}}}}%
    {\raisebox{.14ex}{\vphantom{\ensuremath{\scriptscriptstyle #1}}}}%
    \smash{\hat{#1}}}%
\endgroup%
}
\makeatother

\newcommand{\PP}{\oper{P}}
\newcommand{\QQ}{\oper{Q}}

\def\titlerunning#1{\gdef\@titlerunning{\uppercase{#1}}}\edef\@titlerunning{}%
\def\authorrunning#1{\gdef\@authorrunning{\uppercase{#1}}}\edef\@authorrunning{}%
\def\firstcollaboration#1{\gdef\@firstcollaboration{#1}}\edef\@firstcollaboration{}%
\def\subclass#1{\gdef\@subclass{%
\parbox[t]{\frontmatter@abstractwidth}{{\small \textbf{2010 Mathematical Subject Classification:}~\texttt{#1} \hfill}}{}
}}
\gdef\@subclass{}

\PassOptionsToPackage{pdfauthor={Vladimir V. Kisil},%
    pdftitle={},%
    pdfsubject={mathematics},%
    backref=page,%
  pdfkeywords={symmetries}
}{hyperref}
\usepackage{hyperref}

\newtheorem{thm}{Theorem}
\newtheorem{prop}[thm]{Proposition}
\newtheorem*{principle}{Umbral Principle}
\newtheorem{lem}[thm]{Lemma}

\theoremstyle{definition}
\newtheorem{defn}[thm]{Definition}
\newtheorem{example}[thm]{Example}

    \theoremstyle{remark}
    \newtheorem{rem}[thm]{Remark}

\providecommand{\half}{{\textstyle\frac{1}{2}}}

   \PassOptionsToPackage{british}{babel}
\providecommand{\norm}[2][\relax]{\left\|#2\right\|\ifx#1\relax\else_{#1}\fi}
\providecommand{\modulus}[2][\relax]{\left| #2 \right|\ifx#1\relax\else_{#1}\fi}
\providecommand{\oper}[1]{\mathcal{#1}}
\providecommand{\algebra}[1]{\ensuremath{\mathfrak{#1}}}

\providecommand{\Space}[3][]{\ifx#2R\ifx#1e \mathbb{C}^{#3} \else
\ifx#1p \mathbb{D}^{#3} \else
\ifx#1h \mathbb{O}^{#3} \else
\ifx#1\sigma \mathbb{A}\!^{#3} \else
\ensuremath{\mathbb{#2}^{#3}_{#1}{}} \fi \fi \fi \fi \else
\ensuremath{\mathbb{#2}^{#3}_{#1}{}} \fi}
\providecommand{\FSpace}[3][]{\ensuremath{\ifx#2l \ell_{#3}^{#1}{}\else
  \mathsf{#2}_{#3}^{#1}{}\fi}}

\providecommand{\uir}[3][0]{\ifcase #1{\rho^{#2}_{#3}}%
\or {\breve{\rho}^{#2}_{#3}}%
\or {\tilde{\rho}^{#2}_{#3}}\fi}
\usepackage{braket}
\providecommand{\scalar}[3][\relax]{\ifx#1*\langle #2, #3%
  \rangle \else{\left\langle #2, #3%
    \right\rangle\ifx#1\relax\else_{#1}\fi}\fi}

\providecommand{\SL}[1][2]{\ensuremath{\mathrm{SL}_{#1}(\Space{R}{})}}

\providecommand{\rmi}{\mathrm{i}}
\providecommand{\rme}{\mathrm{e}}
\providecommand{\rmd}{\mathrm{d}}

\makeatletter
\newcommand{\@hslashslash}{%
  \raisebox{-0.2ex}{%
    \scalebox{1}{%
      \rotatebox[origin=c]{20}{$\mathchar'26$%
      }%
    }%
  }%
}
\newcommand{\@kslash}[2]{%
  {%
   \vphantom{#2}%
   \ooalign{\kern#1em\smash{\@hslashslash}\hidewidth\cr
     $#2$\cr
   }%
   \kern.05em
  }%
}
\newcommand\kslash{\mathchoice
  {\hbox{\@kslash{.01}{k}}}
  {\hbox{\@kslash{.01}{k}}}
  {\hbox{\fontsize{\sf@size}{\sf@size}\selectfont\@kslash{.05}{k}}}
  {\hbox{\fontsize{\ssf@size}{\ssf@size}\selectfont\@kslash{.1}{k}}}
}

\providecommand{\myh}{h}
\providecommand{\map}[1]{\mathsf{#1}}

\providecommand{\Zbl}[1]{Zbl\href{http://www.emis.de:80/cgi-bin/zmen/ZMATH/en/zmathf.html?first=1&maxdocs=3&type=html&an=#1&format=complete}{#1}}

\providecommand{\arXiv}[1]{\href{http://arXiv.org/abs/#1}{arXiv:#1}}

\providecommand{\myeprint}[2]{\href{#1}{\texttt{#2}}}
\begin{document}

\titlerunning{Transmutations and umbral calculus} 
\authorrunning{V.V. Kisil} 

\title{Transmutations from 
 the Covariant Transform\\
  on the Heisenberg Group\\
and an Extended Umbral Principle}


\author{\firstname{Vladimir V.} ~\surname{\href{http://www1.maths.leeds.ac.uk/~kisilv/}{Kisil}}}
\email[E-mail: ]{V.Kisil@leeds.ac.uk}
\thanks{On  leave from the Odessa University.}

\affiliation{School of Mathematics,
University of Leeds,
Leeds LS2\,9JT,
England}

\firstcollaboration{(Submitted by S.M. Sitnik)} 

\received{May 24, 2023} 

\begin{abstract} 
  We discuss several seemingly assorted objects: the umbral calculus, generalised translations and associated transmutations, symbolic calculus of operators. The common framework for them is representations of the Weyl algebra of the Heisenberg group by ladder operators. Transporting various properties between different implementations we review some classic results and new opportunities.
\end{abstract}

\subclass{05A40, 34B30, 35S50, 43A85, 44A35, 47G30} 

\maketitle


\medskip
\centerline{\emph{Dedicated to 100\(^{\text{th}}\) anniversary of Ivan Aleksandrovich Kipriyanov's birth}}
  
\bigskip

\hfill\parbox{.3\textwidth}{\small То и это, и вон то,\\ вместе с этим тоже\\---всё равно всегда выходит\\
  как одно и то же.\\
\centerline{Считалочка}}

\section{Introduction: the Bessel operator}
\label{sec:introduction}

Let us start from the celebrated \emph{singular Bessel differential operator}%
\index{singular Bessel differential operator}%
\index{Bessel! singular differential operator}~\cite[\S1.3]{Kipriyanov97a}, \cite[\S2.2.2]{SitnikShishkina19a}, \cite[\S1.4.1]{ShishkinaSitnik20a} given by:
\begin{equation}
  \label{eq:bessel-operator}
  \oper{B}_\nu = \frac{\rmd^2\ }{\rmd t^2} + \frac{\nu}{t} \frac{\rmd \ }{\rmd t} =  t^{-\nu} \, \frac{\rmd \ }{\rmd t} \, t^\nu \, \frac{\rmd \ }{\rmd t}.
\end{equation}
Due to importance of this operator in theoretical and applied settings there are many different approaches to study \(\oper{B}_\nu\) and corresponding boundary value problems. Main tools are various integral transforms \cite{Kipriyanov97a,SitnikShishkina19a,ShishkinaSitnik20a}  intertwining \(\oper{B}_\nu\) with some other operators, which are more accessible for certain reasons. The special r\^ole of such integral transforms is often indicated by calling them \emph{transmutations}\index{transmutation}~\cite{KravchenkoSitnik20a}. 

For example, the \emph{Poisson operator}%
\index{Poisson operator}%
\index{operator!Poisson} is defined for a summable function \(f(t)\) by
\begin{equation}
  \label{eq:poisson-operator}
  [\oper{P}^\nu f](x)= \frac{2 \Gamma (\frac{\nu+1}{2})}{\sqrt{\pi}\Gamma (\frac{\nu}{2}) \, x^{\nu-1}} \int\limits_0^x(x^2-t^2)^{\frac{\nu}{2}-1}\, f(t)\, \rmd t.
\end{equation}
It transmutes (in other words intertwines) the Bessel operator with the second derivative:
\begin{equation}
  \label{eq:poison-intertwining}
  \oper{P}^\nu \circ \oper{B}_\nu = \frac{\rmd^2\ }{\rmd x^2} \circ \oper{P}^\nu.
\end{equation}
Another example is a Fourier-style decomposition over eigenfunctions of \(\oper{B}_\nu\)---the \emph{Hankel transform}%
\index{Hankel transform}%
\index{transform!Hankel} defined by:
\begin{equation}
  \label{eq:hankel-transform}
  [\oper{H}_\nu f](\lambda ) = \int\limits_0^\infty f(t)\, j_{\frac{\nu-1}{2}}(\sqrt{\lambda} t) \,t^\nu\,\rmd t,
\end{equation}
where \(j_{\frac{\nu-1}{2}}(\sqrt{\lambda} t)\) is the solution of the boundary value problem:
\begin{equation}
  \label{eq:little-Bessel-function}
  \oper{B}_\nu j_{\frac{\nu-1}{2}}(\sqrt{\lambda}  t) = -\lambda  j_{\frac{\nu-1}{2}} (\sqrt{\lambda}  t), \qquad \text{ and }
\quad  j_{\frac{\nu-1}{2}}(0)=1, \quad j'_{\frac{\nu-1}{2}}(0)=0.
\end{equation}
The Hankel transform transmutes \(\oper{B}_\nu\) with the operator of multiplication:
\begin{displaymath}
  \oper{H}_\nu \circ \oper{B}_\nu = (-\lambda I) \circ \oper{P}^\nu.
\end{displaymath}

Furthermore, it is fruitful~\cite{Levitan49a,Levitan51a} to extend the  harmonic analysis analogy with \(\oper{H}_\nu\) and introduce the \emph{generalised translations}%
\index{generalised!translations}%
\index{translations!generalised} \(T_t^y\) such that functions \(j_{\frac{\nu-1}{2}}(\sqrt{\lambda} t)\) plays the r\^ole of characters---just like exponents for ordinary translations:
\begin{equation}
  \label{eq:generalised-translation-defn}
  T_t^y j_{\nu}(\lambda t) = j_{\nu}(\lambda y) \, j_{\nu}(\lambda t) .
\end{equation}
There are numerous structural similarities between ordinary and generalised translations~\cite{Levitan49a,Levitan73a}. The correspondence between two resembles the \emph{umbral calculus}%
\index{umbral calculus}%
\index{calculus!umbral}~\cite{Rota75,RotaWay}, that is the technique of dealing with combinatorial quantities \(p_n\) as ``shadows'' of power monomials \(x^n\).  This intriguing magic may be turn into a solid theory through a systematic usage of certain linear operators, see~\cite{Rota75,RotaWay} and \S\ref{sec:co-contr-transf} below.

In this paper the umbral calculus is unfolded through representations of the Weyl algebra by ladder operators. Covariant transforms intertwine properties of corresponding objects from various representations. It allows us to propagate results and methods from one setting to another. We extend the scope through umbral interpretation of generalised shifts and transmutations.  Such connections highlight some gaps and missed opportunities which usually remain in blind spots if the topics are treated in the traditional isolated manner.

The paper \emph{outline} is as follows. We are collecting required preliminary material in Section~\ref{sec:preliminaries}. 
The central place is occupied by representations of Weyl algebra of the Heisenberg group by ladder operators, which are introduced in \S\ref{sec:abstr-repr-weyl}. Then, we remind basics on the covariant (Berezin) transform and Perelomov's coherent states in \S\ref{sec:covar-transf}. The intertwining properties of the covariant transform are foundations of the \emph{umbral principle} presented in~\S\ref{sec:co-contr-transf}.  Convenient complements are the Fourier transform to the momentum picture \S\ref{sec:fourier-transform} and the adjoint action of the ladder operators \S\ref{sec:adjo-acti-resol}. The final bit is a special case of shift invariant delta operator and the resulting binomial in \S\ref{sec:shift-invar-binom}.

Sufficiently equipped by the previous preparations we interpret several examples through ladder operators in Section~\ref{sec:vari-impl}. The combinatorial umbral calculus is our first example in \S\ref{sec:finite-oper-umbr}. We extend it to  the Delsarte--Levitan’s generalised translations and associated transmutations in \S\ref{sec:delsarte-levitan}. These are illustrated by the cases of higher order delta operators---the second derivative in  \S\ref{sec:extend-umbr-calc} and  the Bessel operator \(\oper{B}_\nu\)~\eqref{eq:bessel-operator} in \S\ref{sec:bessel-umbral}. Finally, we discuss an application of symbolic calculus of operators which stems out from the calculus of pseudo-differential operators (PDO) in \S\ref{sec:pdo-type-calculus}.

The final Section~\ref{sec:dicuss-clos-remarks} summarises our finding and proposes an \emph{extended interpretation} of the umbral principle, see illustration~\eqref{eq:umbral-hub}. We also indicate further developments and opportunities generated by the described connections.

\section{Preliminaries}
\label{sec:preliminaries}

Representation theory, in particular of the Heisenberg group~\cite{Howe80a}, is behind many important calculations in
analysis~\cite{Kisil95i,Kisil02a,Kisil12b} and we provide some of its basics in the present section. The
group-theoretical foundations of coherent states/wavelets are
well-known and widely
appreciated~\cite{AliAntGaz00,Perelomov86}. We can widen the applicability of the approach~\cite{Kisil11c} if it is extended to the Banach spaces \cite{FeichtingerGrochenig88a,Kisil98a}.

\subsection{Abstract representation of the Weyl algebra and the Heisenberg group}
\label{sec:abstr-repr-weyl}

 Let an infinite-dimensional vector space \(\FSpace{E}{}\) have a basis \(\{p_n\}\), \(n=0,1,2,\ldots\).  Then, one can define the following associated objects:
 \begin{enumerate}
 \item  The linear map \(\Lower: \FSpace{E}{} \rightarrow \FSpace{E}{}\)  such that:
   \begin{equation}
     \label{eq:hat-Q}
   \Lower p_n = p_{n-1}, \  \text{ for } n>0 \qquad \text { and} \qquad
   \Lower p_0 = 0.
  \end{equation}
  In a combinatorial context \(\Lower\) is called \emph{delta operator}%
  \index{delta operator} or simply \emph{delta}. 
\item The liner map \(\Raise: \FSpace{E}{} \rightarrow \FSpace{E}{} \) such that:
  \begin{equation}
    \label{eq:hat-P}
    \Raise p_n = \iota (n+1)p_{n+1}, \quad \text{ for } n=0,1,2,\ldots,
  \end{equation}
  for some scalar \(\iota\). In many algebraic considerations, e.g. the umbral calculus, it is common to take the simplest value \(\iota=1\). For unitarity in a Hilbert space can employs the imaginary unit \(\iota=\rmi\).   
\item The linear functional \(l_0\in \FSpace[*]{E}{}\) defined by:
  \begin{equation}
    \label{eq:functional-h}
    \scalar[*]{l_0}{ p_n} = \delta_{0n} \,, \qquad \text{ for all } n=0,1,2,\ldots
  \end{equation}
  where \(\delta_{ij}\) is the Kronecker delta. 
\end{enumerate}


The above operators \(\Lower\) and \(\Raise\) satisfy to the \emph{Heisenberg commutation relation}%
\index{Heisenberg!commutation relation}%
\index{commutation relation!Heisenberg} \([\Raise,\Lower]=-\iota I\). Thus, the set \(\{\Raise, \Lower, -\iota I\}\) form a representation of the Weyl Lie algebra \(\algebra{h}\) of the Heisenberg group  \(\Space{H}{}\). In the quantum mechanical framework the action of \(\Raise\)~\eqref{eq:hat-Q} and \(\Lower\)~\eqref{eq:hat-P}  are known as \emph{ladder operators}%
\index{ladder operators}%
\index{operator!ladder}---\emph{creation}%
\index{creation operators}%
\index{operator!creation} and \emph{annihilation}%
\index{annihilation operators}%
\index{operator!annihilation}, respectively~\cite{Folland89,Feinsilver88a,Kisil10a}. Their action is visualised as follows:
\begin{equation}
  \label{eq:ladder-action}
  \xymatrix{ 0 
    &*+[o][F]{\ p_0\ } \ar@/^/[r]|-{\Raise} \ar[l]^{\Lower}  
    &*+[o][F]{\ p_1\  }\ar@/^/[r]|-{\Raise} \ar@/^/[l]|-{\Lower} 
    &*+[o][F]{\ p_2\ }\ar@/^/[r]|-{\Raise} \ar@/^/[l]|-{\Lower}
    &*+[o][F]{\ p_3\ } \ar@/^/[r]|-{\Raise} \ar@/^/[l]|-{\Lower} 
    &*+[o][F.]{\cdots}\ar@/^/@{..>}[l]|-{\Lower}     \ar@/^/@{..>}[r]
     &\ \cdots\ 
    \ar@/^/@{..>}[l]
  }
\end{equation}

We can use the standard exponentiation of  \eqref{eq:hat-Q}--\eqref{eq:hat-P} to obtain an action of the Heisenberg group. Specifically:
\begin{align}
  \label{eq:exp-hat-P}
  \rme^{y \Lower} p_n &= \sum_{k=0}^\infty \frac{y^k \Lower^k}{k!}\, p_n
                            = \sum_{k=0}^n \frac{y^k }{k!}\, p_{n-k} ; \\
  \label{eq:exp-hat-Q}
    \rme^{x \Raise} p_n &= \sum_{k=0}^\infty \frac{x^k \Raise^k}{k!}\, p_n
                              = \sum_{k=0}^\infty \binom{n+k}{k} x^k\,  p_{n+k} .
\end{align}
One may notice a difference between~\eqref{eq:exp-hat-P} and~\eqref{eq:exp-hat-Q}: the former is effectively a finite linear combination of \(p_n\) and the latter is infinite. Thus, \(\rme^{y \Lower}\) makes sense in the linear space \(\FSpace{E}{}\) itself, but \(\rme^{x \Raise}\) needs some additional structure, e.g. a topology to interpret~\eqref{eq:exp-hat-Q} as a convergent series. An undemanding resolution is the ``formal power series'' framework, that is a collection of all infinite series \(f=\sum_0^\infty a_n p_n\) with \(f^{(k)} \rightarrow f\) if and only if \(a_n^{(k)} \rightarrow a_n\) for all \(n\). We will continue to loosely denote some topological extension of \(\FSpace{E}{}\) which supports expansions~\eqref{eq:exp-hat-P}--\eqref{eq:exp-hat-Q} by the same letter \(\FSpace{E}{}\).

Once actions~\eqref{eq:exp-hat-P}--\eqref{eq:exp-hat-Q} are meaningful, they satisfy to  the \emph{Weyl commutator relation}%
\index{Weyl!commutation relation}%
\index{commutation relation!Weyl}: 
\begin{displaymath}
\rme^{y \Lower}\, \rme^{x \Raise} = \rme^{\iota xy}\, \rme^{x \Raise}\,\rme^{y \Lower}.
\end{displaymath}
Thus, we have a representation of the (polarised) \emph{Heisenberg group}%
\index{Heisenberg!group}%
\index{group!Heisenberg} \(\Space{H}{}\) of the form:
\begin{equation}
  \label{eq:H-representation}
  \uir{}{\FSpace{E}{}}(s,x,y) =  \rme^{-\iota s}\,\rme^{-y \Lower}\, \rme^{-x \Raise},
\end{equation}
with the composition law reflecting the multiplication on \(\Space{H}{}\):
\begin{equation}
  \label{eq:H-group-law}
  \uir{}{\FSpace{E}{}}(s,x,y) \uir{}{\FSpace{E}{}}(s',x',y') = \uir{}{\FSpace{E}{}}(s+s'+xy',x+x',y+y'),
\end{equation}
\begin{example}
  \label{ex:archetypal-polinom}
  The \emph{archetypal model}%
  \index{archetypal model}%
  \index{model!archetypal }
  is the linear space of polynomials \(\FSpace{\underline{E}}{}\), say, in a variable \(t\) with the monomial basis \(\underline{p}_n(t) = \frac{1}{n!}t^n\). Then \(\underline{\Lower}= \frac{\rmd\ }{\rmd t}\), \(\underline{\Raise}= \iota t I \) (the operator of multiplication by the variable) and \(\scalar[*]{\underline{l}_0}{p} = p(0)\) (the evaluation at \(0\)). An extension of \(\FSpace{\underline{E}}{}\) to a space of power series in \(t\) with rapidly decreasing coefficients allows us to write the representation~\eqref{eq:H-representation} as:
  \begin{equation}
    \label{eq:schrodinger-repres}
    \uir{}{\iota}(s,x,y) f(t) = \rme^{-\iota(s+x(t-y))} f(t-y),
  \end{equation}
  which is a cousin of the Schr\"odinger representation~\cite[\S1.3]{Folland89}. 
\end{example}

\subsection{Covariant transform and intertwining properties}
\label{sec:covar-transf}
There is an extensive literature on the covariant transform (also known as the coherent state transform and many other names) in general~\cite[Thm.~8.1.3]{AliAntGaz00} and for the Heisenberg group in particular~\cite{Folland89,Kisil10a,AlameerKisil21a}. To make a long story short we recall here only the bare minimum of definitions and notations, see provided references for further details.

It is common to consider covariant transform either in inner product spaces or Gelfand triples  (ridged Hilbert spaces)~\cite{FeichtingerGrochenig88a}. In the present situation it is more convenient to work in a linear space and its dual~\cite{Kisil98a,Kisil12d}. 
\begin{defn} \cite{Kisil09d} Let \(\uir{}{}\) be a representation of
  a group \(G\) in a vector space \(\FSpace{V}{}\). Take a   \emph{fiducial functional}%
  \index{fiducial functional}%
  \index{functional!fiducial}
 \(l\in  \FSpace[*]{V}{}\) and denote its pairing with \(v \in \FSpace{V}{}\) by \(\scalar{l}{v}\).  The \emph{covariant transform}%
  \index{covariant!transform}%
  \index{covariant!transform}%
  \index{transform!covariant} is the map:
  \begin{equation}
    \label{eq:covariant-trans}
    [\oper{W}_l v](g)=\scalar{l}{\uir{}{}(g)v},\qquad v\in \FSpace{V}{},
    \quad g\in   G,
  \end{equation}
  to scalar-valued functions on \(G\).
  In the case of a Hilbert
space \(\FSpace{V}{}\), a fiducial functional is provided by a pairing with a vector
\(f\in \FSpace{V}{} \sim \FSpace[*]{V}{}\) known as a \emph{mother wavelet}%
\index{wavelet!mother}%
\index{mother wavelet} or \emph{vacuum state}%
\index{vacuum state}%
\index{state!vacuum}~\cite{AliAntGaz00,Perelomov86}.
\end{defn}


The covariant transform \(\oper{W}_l\) plainly interacts with the left \(\Lambda\)%
\index{left regular representations!}%
\index{regular!representations!left}%
\index{representations!regular!left} and right \(R\) \emph{regular representations}%
\index{regular!representation!right}%
\index{representation!regular!right}%
\index{right!regular representation} of \(G\), which are:
\begin{displaymath}
  \Lambda(g):   f(h)\mapsto f(g^{-1}h) \qquad \text { and } \qquad
  R(g): f(h)\mapsto f(hg), \qquad \text{ where } h,g \in G.
\end{displaymath}

\begin{lem}
  \label{le:covariant-intertwining}
  The covariant transform~\eqref{eq:covariant-trans} intertwines the left \(\Lambda\) and right \(R\) regular
  representations of \(G\) with the following
  actions:
  \begin{equation}
    \label{eq:covariant-intertwine}
    \Lambda  (g) \circ \oper{W}_l = \oper{W}_l \circ \uir{}{}(g)
    \quad\text{ and }\quad
     R (g) \circ \oper{W}_l= \oper{W}_{ \uir{*}{}(g) l }
    \quad\text{ for all } g\in G.
  \end{equation}
  Here \(\uir{*}{}\) is the adjoint representation on \(\FSpace[*]{V}{}\), that is \(\scalar[*]{l}{\uir{}{}(g) v} = \scalar[*]{\uir{*}{}(g)l}{ v}\).
\end{lem}


For a subgroup \(H\subset G\) and its character \(\chi\) let the fiducial functional \(l\) has the \emph{covariance property}:
\begin{equation}
  \label{eq:covariance-functional}
  \scalar{l}{\uir{}{}(h) v} = \bar{\chi}(h)\scalar{l}{ v},\qquad \text{ for all } h\in H
  \quad \text{ and all } v \in \FSpace{V}{}. 
\end{equation}
Then, \(\oper{W}_l\) is a Perelomov-style covariant transform~\cite{Perelomov86}, that is only a part of the covariant transform~\eqref{eq:covariant-trans} allows to recover all values.  To this end, we fix a continuous section \(\map{s}:
G/H \rightarrow G\), which is a right inverse to the
natural projection \(\map{p}: G \rightarrow G/H\).
\begin{defn} 
  \cite[\S5.1]{Kisil11c} Let \(l \in \FSpace[*]{V}{}\) intertwine the restriction of
  \(\uir{}{}\) to \(H\)  with a character \(\chi\) of \(H\), cf.~\eqref{eq:covariance-functional}. Then, the
  \emph{induced covariant transform}%
  \index{induced!covariant transform}%
  \index{covariant!transform!induced}%
  \index{transform!covariant!induced} is:
  \begin{equation}
    \label{eq:induced-covariant-trans}
    [\oper{W}_l v](x)=\scalar{l}{\uir{}{}(\map{s}(x))v}, \qquad v\in \FSpace{V}{},
    \quad x\in    G/H.
  \end{equation}
\end{defn}
Under our assumptions, the induced covariant transform intertwines
\(\uir{}{}\) with a \emph{representation induced}%
\index{representation!induced}%
\index{induced!representation} from \(H\) by the
character \(\chi\). 

In many cases, e.g. for square integrable representations%
\index{square integrable!representations}%
\index{representations!square integrable} and an admissible mother
wavelet \(v\in \FSpace{V}{}\), the image space of the covariant transform is a
reproducing kernel Hilbert space~\cite[Thm.~8.1.3]{AliAntGaz00}. That
means that for any function \(l \in \oper{W}_l \FSpace{V}{}\) we have the
integral reproducing formula:
\begin{equation}
  \label{eq:reproducing-integral}
  v(y)=\int_X v(x)\,\bar{k}_y(x)\,dx,
\end{equation}
where the \emph{reproducing kernel}%
\index{reproducing kernel}%
\index{kernel!reproducing} \(k_y\) is the twisted convolution
with the normalised covariant transform \(\oper{W}_l
(\uir{}{}(\map{s}(y)^{-1}) l)\) for the mother wavelet \(l\)~\cite{Kisil13a}. 

\subsection{Co
  variant transforms for the Heisenberg group and the Umbral principle }
\label{sec:co-contr-transf}

We can observe that for the linear functional \(l_0\)~\eqref{eq:functional-h} the identity holds:
\begin{displaymath}
  \scalar[*]{l_0}{\rme^{x \Raise} f} = \rme^{x} \scalar[*]{l_0}{ f},
\end{displaymath}
which can be initially verified for \(f=p_n\) and then extended to the whole \(\FSpace{E}{}\) by linearity. Therefore, there is an induced covariant transform~\eqref{eq:induced-covariant-trans} for the subgroup
\begin{equation}
  \label{eq:subgroup-Hy}
  H_y=\{(s,x, 0) \mid s,x \in \Space{R}{}\} \quad \text{ and the respective homogeneous space } \quad   \Space{H}{}/H_y \sim \Space{R}{}.
\end{equation}
The representation \(\uir{}{\myh}\) of \(\Space{H}{}\) induced from the character \(\chi_{\myh}(s,x,0) = \rme^{\iota s}\) of \(H_y\)  coincides with the Schr\"odinger representation~\eqref{eq:schrodinger-repres}~\cite{Kisil10a}:
\begin{displaymath}
  [\uir{}{\myh}(s,x,y) f](u) = \rme^{-\iota(s+x(u-y))} f(u-y).
\end{displaymath}
Then, the induced covariant transform defined by functional \(l_0\)~\eqref{eq:functional-h} and the subgroup \(H_y\)
\eqref{eq:subgroup-Hy} maps \(\FSpace{E}{}\) to a space of functions of one variable: 
\begin{equation}
  \label{eq:covariant-umbral}
 \tilde{f}(u) \coloneqq [\oper{W}_0 f](u)=\scalar[*]{l_0}{\rme^{u \Lower} f}, \qquad \text{ in particular } \quad \tilde {p}_n(u) = \scalar[*]{l_0}{\rme^{u \Lower} p_n} = \frac{(\iota u)^n}{n!},
\end{equation}
where the last identity follows from~\eqref{eq:exp-hat-P}. The transform \(\oper{W}_0\) intertwines the representation \(\uir{}{\FSpace{E}{}}\)~\eqref{eq:H-representation} with the representation \(\uir{}{\iota}\)~\eqref{eq:schrodinger-repres}, in particular:\footnote{The archetypal implementation of a covariant transform and these intertwining relations is  the Fock--Segal--Bargmann transform, cf. ~\cite{Kisil19b}. It intertwines the creation and annihilation operators for a quantum harmonic oscillator with operators of multiplication by the complex variable \(z\) and the complex derivative \(\frac{\rmd\ }{\rmd z}\), respectively.}
\begin{equation}
  \label{eq:intertwining}
  \oper{W}_0 \circ \uir{}{\FSpace{E}{}} = \uir{}{\iota} \circ \oper{W}_0 , \qquad
  \oper{W}_0 \circ \Lower = \frac{\rmd\ }{\rmd u} \circ \oper{W}_0 , \qquad
  \oper{W}_0 \circ \Raise = (\iota uI) \circ \oper{W}_0 , \qquad
  \scalar[*]{l_0}{f} = [\oper{W}_0f](0).
\end{equation}

The above relations are the fundamentals of the following principle, cf. \S\ref{sec:finite-oper-umbr} below:
\begin{principle}
  Any statement on functions of one variable, which is formulated in terms of a linear combination of derivatives, multiplication by monomials and evaluation at \(0\), corresponds to a statement about elements of \(\FSpace{E}{}\) expressed through \(\Lower\), \(\Raise\) and \(l_0\) according to the vocabulary~\eqref{eq:intertwining}.
\end{principle}

Note that the above principle does not require any additional assumptions about the nature of \(\FSpace{E}{}\) or its basis \(\{p_n\}\). Once they are fixed, operators \(\Lower\)~\eqref{eq:hat-Q}, \(\Raise\)~\eqref{eq:hat-P} and the functional \(l_0\)~\eqref{eq:functional-h} are completely defined and the Umbral principle is fully set.

\begin{rem}
  \label{re:role-Heisenberg-group}
  Interestingly, the umbral correspondence \(\oper{W}_0\)~\eqref{eq:intertwining}  is defined only through the exponentiation of the delta \(\Lower\) in~\eqref{eq:covariant-umbral}. Thus, the paired operator \(\Raise\) is missing or at least obscured in many considerations. Therefore, it is tempting to deem the operator \(\Raise\) to be optional or even excessive. However, the importance of \(\Raise\) manifests itself through the possibility to define the entire umbral framework \((\FSpace{E}{}, \{p_n\}_0^\infty)\) through the triple \((\FSpace{E}{}, \Raise, p_0)\). Indeed, the sequence \(\{p_n\}_0^\infty\) is inductively defined, cf.~\eqref{eq:exp-hat-Q} by \(p_{n+1}= \frac{1}{\iota (n+1)} \Raise p_n\), \(n=0, 1, 2,\ldots\).
  On the other hand, there is no way equally well describe the situation through the delta \(\Lower\) alone, see the example of different Appell polynomials, which all share the same \(\Lower=\frac{\rmd\ }{\rmd t}\), in \S\ref{sec:finite-oper-umbr} below. We will meet some more arguments later, e.g. a usage of the operator \(\Raise\) allows an extension of the operator calculus in \S\ref{sec:pdo-type-calculus}.
\end{rem}

\subsection{The momentum picture and the Fourier transform}
\label{sec:fourier-transform}

In some circumstances it is preferable to have an intertwining property with operators  \(\Lower\) and \(\Raise\)  swapped in comparison to~\eqref{eq:intertwining}. In the physical language: to use momenta of a particle instead of its coordinates in the configuration space. To this end we replace the functional \(l_0\)~\eqref{eq:functional-h} by a functional \(m\) invariant under the operator \(\rme^{u \Lower}\) 
\begin{equation}
  \label{eq:functional-inv}
  \scalar[*]{m}{ \rme^{u \Lower} f} = \scalar[*]{m}{f}\,.
\end{equation}
Then, the respective induced covariant transform:
\begin{equation}
  \label{eq:intertwining-second}
  [\oper{W}_{m} f](v) = \scalar[*]{m}{\rme^{v \Raise} f} , \ 
  \text{ intertwines } \ 
  \oper{W}_{m} \circ \Lower = (-\iota vI) \circ \oper{W}_{m} , \text{ and }
  \oper{W}_{m} \circ \Raise = \frac{\rmd\ }{\rmd v} \circ \oper{W}_{m} .
\end{equation}

The above consideration is purely formal because the functional~\eqref{eq:functional-inv} does not have a non-trivial bounded action on the basis \(p_n\) of \(\FSpace{E}{}\). Indeed,  \(\rme^{u \Lower} p_1 = p_1+ u p_0\) and we shall have \(\scalar[*]{m}{p_1} = \scalar[*]{m}{p_1} +u \scalar[*]{m}{p_0}\) for all \(u\). Boundedness on \(p_n\) implies \(\scalar[*]{m}{p_0}=0\), which can be extended to \(\scalar[*]{m}{p_n}=0\) for all \(n\) by induction.

Yet, such a non-zero functional may exist on a certain subspace of the extended space \(\FSpace{E}{}\). For example, if \(\Lower=\frac{\rmd \ }{\rmd t}\) and  \(\rme^{u \Lower}: f(t) \mapsto f(t+u)\), cf. Ex.~\ref{ex:archetypal-polinom}, a  shift-invariant functional is:
\begin{displaymath}
  \scalar[*]{m}{f}= \int\limits_{-\infty}^\infty f(t)\, \rmd t \qquad
  \text{ for } f \in\FSpace{L}{1}(\Space{R}{}).
\end{displaymath}
Then, the induced covariant transform \(\oper{W}_{m}\)~\eqref{eq:intertwining-second} effectively becomes the Fourier transform:
\begin{displaymath}
  [\oper{W}_{m}f](v) = \int\limits_{-\infty}^\infty \rme^{-\iota vt}  f(t)\, \rmd t,
\end{displaymath}
It intertwines the Schr\"odinger representation \(\uir{}{\iota}\)~\eqref{eq:schrodinger-repres} with itself through  
an automorphism of \(\Space{H}{}\)~\cite{Howe80a}, \cite[\S1.3]{Folland89}. On the level of the Weyl algebra representation spanned by operators \(\Lower= \frac{\rmd\ }{\rmd t}\), \(\Raise= \iota t I \) from Ex.~\ref{ex:archetypal-polinom} the automorphism acts as follows:
\begin{displaymath}
  \oper{W}_{m}: \Lower \mapsto -\Raise, \quad \text{ and } \quad
  \oper{W}_{m}: \Raise \mapsto \Lower\,, \qquad \text{ therefore } \qquad
  [\oper{W}_{m} \Raise, \oper{W}_{m} \Lower] = [\Raise, \Lower] =-\iota I.
\end{displaymath}
The above discussed incompatibility of the averaging functional  \(m\) and the basis \(\{p_n\}\) is just another reason why it is preferable to set the umbral framework through operators \(\Lower\) and \(\Raise\) rather than trough a specific basis \(\{p_n\}\), cf. Rem.~\ref{re:role-Heisenberg-group}.

\subsection{Adjoint action, resolution of the identity and binomial formula}
\label{sec:adjo-acti-resol}

We define adjoint ladder operators \(\FSpace[*]{E}{} \rightarrow \FSpace[*]{E}{} \) in the usual way:
\begin{equation}
  \label{eq:adjoint-op-defn}
  \scalar[*]{\Lower^*l}{f} = \scalar[*]{l}{\Lower f}, \quad
  \scalar[*]{\Raise^*l}{f} = \scalar[*]{l}{\Raise f} \quad \text{ for all }
  f\in\FSpace{E}{} \text{ and } l\in \FSpace[*]{E}{}.
\end{equation}
Then, the sequence \(\{l_k\}_0^\infty\subset \FSpace[*]{E}{}\) is inductively produced  starting from the functional \(l_0\)~\eqref{eq:functional-h}:
\begin{equation}
  \label{eq:functional-sequence}
  l_{n+1} = \Lower^* l_n, \qquad n=0,1,2,\ldots.
\end{equation}
Note, that the passage to adjoint operators swaps the creation and annihilation r\^oles of ladder operators. 

The properties~\eqref{eq:functional-h} and~\eqref{eq:adjoint-op-defn} implies bi-orthogonality of sequences \(\{p_n\}_0^\infty\) and \(\{l_k\}_0^\infty\):
\begin{equation}
  \label{eq:bi-orthogonal}
  \scalar[*]{l_k}{p_n} = \delta_{kn}.
\end{equation}
Therefore, we have the following \emph{resolution of the identity}%
\index{resolution of the identity}%
\index{identity, resolution of} written in the Dirac bra-ket notation:
\begin{equation}
  \label{eq:resolution-identity}
  I=\sum_{k=0}^\infty \ket{p_k}\bra{l_k}, \qquad \text{ that is }
  \qquad
  f = \sum_{k=0}^\infty \scalar{l_k}{f}\, p_k, \qquad \text{ for all } f \in \FSpace{E}{}.
\end{equation}
The last identity for \(f=p_n\) directly follows from~\eqref{eq:bi-orthogonal} and then extends to all elements by linearity.

One can attempt an \emph{unitarisation trick}: define a map \(\oper{F}: \FSpace{E}{} \rightarrow \FSpace[*]{E}{}\) by the rule \(\oper{F}: p_n \mapsto l_n\) for all \(n\). Say, if \(\FSpace{E}{}\) is a space of polynomials in a single variable, we may look for a measure \(\mu\), which implements~\eqref{eq:bi-orthogonal} in the form:
\begin{equation}
  \label{eq:orthogonal-polynomials}
  \scalar{\oper{F} p_k}{p_n} = \int p_k(t)\, p_n(t)\,\rmd\mu(t) = \delta_{kn}.
\end{equation}
That is the classical approach to the \emph{orthogonal polynomials}%
\index{orthogonal polynomials}%
\index{polynomial!orthogonal} from ladder operators~\cite{Feinsilver88a}.
Let we have two abstract representations \((\FSpace{E}{}, \{p_n\}_0^\infty, \Lower, \Raise, l_0)\) and \((\tilde{\FSpace{E}{}}, \{\tilde{p}_n\}_0^\infty, \tilde{\Lower}, \tilde{\Raise}, \tilde{l}_0)\). Then we may define the correspondence \(\FSpace{E}{} \rightarrow \tilde{\FSpace{E}{}}\) 
similarly to the resolution of identity~\eqref{eq:resolution-identity}: 
\begin{equation}
  \label{eq:Intertwining-E-E1}
    \oper{V}=\sum_{k=0}^\infty \ket{\tilde{p}_k}\bra{l_k}, \qquad \text{ that is }
  \qquad
  f = \sum_{k=0}^\infty \scalar{l_k}{f}\, \tilde{p}_k, \qquad \text{ for all } f \in \FSpace{E}{}.
\end{equation}
Clearly, \eqref{eq:functional-sequence} implies \(\oper{V}: p_k \mapsto \tilde{p}_k\) for all \(k=0,1,2,\ldots\).

Similarly we can define \emph{generating function}%
\index{generating function}%
\index{function|generating} between two representations \((\FSpace{E}{}, \{p_n\}_0^\infty)\) and \((\tilde{\FSpace{E}{}}, \{\tilde{p}_n\}_0^\infty)\):
\begin{equation}
  \label{eq:generating-function}
  F(s,t) = \sum_{k=0}^\infty k! \, \tilde{p}_k(s) \, p_k(t), \qquad \text{ which intertwines } \quad  \Lower_t F(s,t) =\tilde{\Raise}_s F(s,t).
\end{equation}
In most cases the generating function is taken for \((\tilde{\FSpace{E}{}}, \{\tilde{p}_n\}_0^\infty)\) being the archetypal model \((\FSpace{\underline{E}}{}, \{\underline{p}_n\}_0^\infty)\) from Ex.~\ref{ex:archetypal-polinom}, that is:
\begin{equation}
  \label{eq:generating-function-archetypal}
  F(s,t) = \sum_{k=0}^\infty s^k p_k(t), \qquad \text{ such that } \quad \Lower_t F(s,t) = s F(s,t).
\end{equation}

\subsection{Shift invariance, binomial formula and operator calculus}
\label{sec:shift-invar-binom}

There is a  particular but still important situation of  a translation-invariant delta \(\Lower\).
\begin{prop}
  \label{pr:binomial-type}
  Let \(\FSpace{E}{}\) be a space of function on the real line, let \(\Lower\) commute with all translations \(T^y: f(t) \mapsto f(t+y)\), \(y\in\Space{R}{}\) and \(\scalar[*]{l_0}{f}=f(0)\). Then
\begin{equation}
  \label{eq:binomial-abstract}
  p_n(t+y) = \sum_{k=0}^n p_{n-k}(y)\, p_k(t).
\end{equation}
Therefore, we can express the translation \( T^y=\rme^{y \frac{\rmd \ }{\rmd t}}\) as a function of \(\Lower\):
\begin{equation}
  \label{eq:shift-through-binomial}
  T^y f= \sum_{k=0}^{\infty}  p_k(y) \, \Lower^k f.
\qquad \text{ for any } f \in \FSpace{E}{}.
\end{equation}
\end{prop}
A proof begins from an application of~\eqref{eq:resolution-identity} to \(T^y p_n(t)=p_n(t+y)\): 
\begin{displaymath}
  \begin{split}
  p_n(t+y) &= \sum_{k=0}^\infty \scalar{l_k}{T^y p_n}\, p_k(t) = \sum_{k=0}^\infty \scalar[*]{\Lower^{*k}l_0}{ T^y p_n}\, p_k(t)= \sum_{k=0}^\infty \scalar[*]{l_0}{\Lower^k T^y p_n}\, p_k(t)\\
  &= \sum_{k=0}^\infty \scalar[*]{l_0}{T^y \Lower^k p_n}\, p_k(t) = \sum_{k=0}^n \scalar{l_0}{T^y p_{n-k}}\, p_k(t)= \sum_{k=0}^n p_{n-k}(y)\, p_k(t),
\end{split}
\end{displaymath}
where the series obviously terminates at \(k=n\). Then, \eqref{eq:shift-through-binomial} for \(f=p_n\) is exactly~\eqref{eq:binomial-abstract}, thereafter \eqref{eq:shift-through-binomial} for a general \(f \in \FSpace{E}{}\) follows from linearity.

In the archetypal case of Ex.~\ref{ex:archetypal-polinom} with \(p_n(t) = \frac{1}{n!}t^n\) identity~\eqref{eq:binomial-abstract} turns out to be the celebrated binomial formula:
\begin{equation}
  \label{eq:binomial-formula}
  \frac{(t+x)^n}{n!} = \sum_{k=0}^n  \frac{x^{n-k}}{(n-k)!}\,  \frac{t^{k} }{k!}  
                                                     \qquad \text{ or } \qquad
  (t+x)^n= \sum_{k=0}^n{ n \choose k} x^{n-k}\,  t^k.
\end{equation}
Yet, \eqref{eq:binomial-abstract} remains valid for numerous polynomials of binomial types, see below.


Now let us turn to some examples of the above abstract scheme.

\section{Various implementations}
\label{sec:vari-impl}
We start from the original version of the umbral calculus and then extend the scope by some new illustrations. 

\subsection{Finite Operator (Umbral) Calculus}
\label{sec:finite-oper-umbr}

As it often happens, the umbral calculus started from some particular observations in specific circumstances. Thereafter, it took several iterations to separate an abstract core from technical aspects~\cite{Rota95} and recognise the r\^ole of the Heisenberg group~\cite{Cigler78,Feinsilver88a,Kisil97b,Kisil01b,Kwasniewski03a,LeviTempestaWinternitz04a,DattoliLeviWinternitz08a,FaustinoRen13a} in the construction. Here we are moving the opposite direction: from the abstract scheme of \S\ref{sec:abstr-repr-weyl} to its specific implementations.

It is common in combinatorics to take \(p_n\) to be a polynomial of degree \(n\). Alternatively, it can be requested that the delta operator \(\Lower\) sends the first order monomial to a constant: \(\Lower t = c\)~\cite{FoundationIII}. Thereafter, some  additional assumptions are employed as well, e.g. in the reverse historical order:
\begin{itemize}
\item If a delta \(\Lower\) commutes with shifts then \(\{p_n\}\) is called a \emph{Sheffer sequence}.%
  \index{Sheffer polynomials}%
  \index{polynomial!Sheffer}
\item If \(\Lower=\frac{\rmd\ }{\rmd t}\) then  \(\{p_n\}\) is a sequence of \emph{Appell polynomials}%
  \index{Appell polynomials}%
  \index{polynomials!Appell} (thus, they are special case of Sheffer polynomials)
\item If Sheffer polynomials satisfy \(\scalar[*]{l_0}{p_n} = p_n(0)\) then \(\{p_n\}\)  are called sequences of \emph{binomial type}%
  \index{polynomials!binomial type}%
  \index{binomial type}. As we already know, such \(\{p_n\}\) shall satisfy to the binomial formula~\eqref{eq:binomial-abstract}, which can be taken as their alternative definition and, clearly, is the source of their name.
\end{itemize}
Notably, the only Appell polynomials of binomial type are monomials \(p_n(t) = \frac{1}{n!} (\iota t)^n\) from Ex.~\ref{ex:archetypal-polinom}. That is the Appell and binomial type polynomials are two different branches springing from the archetypal model and still covered by the same umbrella of Sheffer polynomials.  

\begin{example}
There are numerous sequences of polynomials covered in the literature, cf.~\cite{Rota95,RotaWay}: 
\noindent
\begin{center}
\begin{tabular}{||l|c|c|c|c||}
  \hline
  \hline
&  $\frac{\strut}{\strut}p_n(t)$ & $ \Lower$ & $\Raise$ &  $\scalar[*]{l_0}{f}$ \\
  \hline
Monomials  &
  $\frac{\strut}{\strut}\frac{1}{n!}t^n$ & $\frac{\rmd\ }{\rmd t}$ & $tI$ &  $f(0)$ \\
  \hline
Lower factorials &
   $\frac{\strut}{\strut}\frac{1}{n!} (t)_n= \frac{1}{n!} t(t-1)(t-2)\cdots(t-n+1)$ & $f(t+1)-f(t)$ & \eqref{eq:hat-P} &  $f(0)$ \\
  \hline
  Upper factorials&  $\frac{\strut}{\strut}\frac{1}{n!} t^{(n)} = \frac{1}{n!} t(t+1)(t+2)\cdots(t+n-1)$ & $f(t)-f(t-1)$ & \eqref{eq:hat-P} &  $f(0)$ \\
  \hline
  \ldots & \ldots &\ldots &\ldots &\ldots \\
  \hline
  \hline
\end{tabular}
\end{center}
with more classical names to follow: the Abel polynomials, the Touchard polynomials, Appell sequences, Hermite polynomials, Bernoulli polynomials, etc. In many cases the simplest description of the operator \(\Raise\) is given by the references to~\eqref{eq:hat-P}. Yet a sort of analytic expressions may be elaborated sometimes as well.
\end{example}

We illustrate the Umbral principle here with one example only. Rewrite the binomial formula~\eqref{eq:binomial-formula} using the operator \(\frac{\rmd\ }{\rmd x}\):
\begin{align*}
  \rme^{(y \frac{\rmd\ }{\rmd x})} \frac{x^n}{n!} = \sum_{k=0}^n   \frac{y^{k} }{k!} \frac{x^{n-k}}{(n-k)!} 
  = \sum_{k=0}^n   \frac{y^{k} }{k!} \left(\frac{\rmd^k}{\rmd x^k}\,  \frac{x^n}{n!}\right). 
\end{align*}
Then an application of the Umbral principle (e.g. the covariant transform) for Appell polynomials \(p_n\)  with \(\Lower=\frac{\rmd\ }{\rmd x}\) gives the identity:
\begin{displaymath}
  p_n(t+y) =  \rme^{y\Lower} p_n(t) =  \sum_{k=0}^n   \frac{y^{k} }{k!} \, \Lower^k  p_n(t)
  =  \sum_{k=0}^n   \frac{y^{k} }{k!} \,  p_{n-k}(t). 
\end{displaymath}
More illustrations 
can be found in~\cite{Rota95,RotaWay}.

\subsection{Delsarte--Levitan's generalised translations}
\label{sec:delsarte-levitan}

Now we link the umbral principle with
\emph{generalised translations}%
\index{generalised translation}%
\index{translation!generalised} proposed by Delsarte and extensively investigated by Levitan, cf. \cite{Levitan73a}, \cite[\S2.2]{SitnikShishkina19a}, \cite[\S3.4.3]{ShishkinaSitnik20a}. We present an adaptation of the original Delsarte's approach, which is well tuned to our theme.

Let \(L\) be a linear operator on a space \(\FSpace{E}{}(\Space{R}{})\) of functions in one variable. We also assume that for some neighbourhood \(\Omega\) of \(0\) and for any scalar \(\lambda \in \Omega\) there is a unique non-zero solution \(\phi(\lambda, t)\) of the eigenvalue problem:
\begin{equation}
  \label{eq:L-eigenvalue-problem}
  \begin{cases}
    L_t\phi(\lambda, t) = \lambda\, \phi(\lambda , t)\quad \text{ (where \(L_t\) acts in variable \(t\))};\\
    \phi(\lambda, 0) = 1 .
  \end{cases}
\end{equation}
Furthermore, let \(\phi(\lambda, t)\) be analytic at \(\lambda=0\) and generates functions \(\phi_n(t)\) in the power expansion:
\begin{equation}
  \label{eq:generating-phi-n}
  \phi(\lambda, t) = \sum_{n=0}^\infty \lambda ^n \, \phi_n(t) \quad \text{  such that \eqref{eq:L-eigenvalue-problem} implies } \quad L\phi_n =\phi_{n-1} \quad \text{ and } \quad \phi_n(0)=\delta_{0n} 
  .
\end{equation}
In other words, the collection \((\FSpace{E}{}(\Space{R}{}), \{\phi_n\}_0^\infty)\) provides the abstract representation of the Heisenberg group from \S\ref{sec:abstr-repr-weyl} with \(\Lower=L\) and \(\scalar{l_0}{f}=f(0)\). Furthermore, in this setting the decomposition~\eqref{eq:generating-phi-n} is an implementation of the generating function~\eqref{eq:generating-function-archetypal}.

Changing the boundary conditions \(\phi(\lambda, 0) = 1\) in~\eqref{eq:generating-phi-n} we get different umbral sequences \(\{\phi_n\}_0^\infty\) for the same operator \(L\). This corresponds to various Appell polynomials for the same \(\Lower = \frac{\rmd\ }{\rmd t}\) in the classic umbral setting outlined above.

If \(L\) commutes with all ordinary translations \(\rme^{y\frac{\rmd\ }{\rmd t}}: f(t) \mapsto f(t+y)\), then \(\rme^{y\frac{\rmd\ }{\rmd t}} = \sum_{0}^{\infty}  \phi_k(y) L_t^k\) by~\eqref{eq:shift-through-binomial}. If translation-invariance of \(L\) is not assumed then the previous formula defines the umbral version for the \emph{generalised translation}%
\index{generalised translation}%
\index{translation!generalised} \(T_t^y\):
\begin{equation}
  \label{eq:generalised-translation}
  [T_t^y f](t)= \sum_{0}^{\infty}  \phi_k(y)\, L^k f(t).
\end{equation}
By induction the intertwining property~\eqref{eq:generating-function-archetypal} of the generating function \(\phi(\lambda,t)=\sum_0^\infty \lambda^k \phi_k(t)\)~\eqref{eq:generating-phi-n} implies \(L^k \phi(\lambda,t) = \lambda^k  \phi(\lambda,t)\), therefore, cf.~\eqref{eq:generalised-translation-defn}:
\begin{equation}
  \label{eq:generalised-character-property}
  [T_t^y \phi](\lambda,t) =\sum_{k=0}^{\infty}  \phi_k(y)\, L^k \phi(\lambda,t) = \sum_{k=0}^{\infty}  \phi_k(y)\, \lambda^k  \phi(\lambda,t)
  =\phi(\lambda,y)\,\phi(\lambda,t).
\end{equation}
That is, the generating function \(F(s,t)\)~\eqref{eq:generating-function-archetypal} is a character of the generalised translation~\eqref{eq:generalised-translation}.   Now we turn to specific examples of generalised translations.

\subsubsection{Second order derivative}
 \label{sec:extend-umbr-calc}

To begin with, we can drop one of the main assumption of the umbral calculus in combinatorics that \(p_n\) is a polynomial of degree \(n\).
For example, let us consider even-order monomials \(p_n(t) =\frac{(-1)^n\,t^{2n}}{(2n)!}\), \(n=0,1,2,\ldots\) and the linear space \(\FSpace{E}{2}\) spanned by them. The second derivative \(\QQ = 
-\frac{\rmd^2\ }{\rmd t^2} \) has the action \(\QQ p_n =  p_{n-1}\) for \(n>1\) and \(\QQ p_0 = 0\). Therefore, to have the  commutator \([\PP,\QQ]= I\) we can 
define a linear operator \(\PP\) by the rule
\begin{displaymath}
  \PP p_n = (n+1) p_{n+1}, \qquad \text{ that is } \qquad \PP t^{2n} = -\frac{1}{2(2n+1)}t^{2n+2}.
\end{displaymath}
Alternatively, we can employ antiderivatives and get an analytic expression, cf.~\cite{DattoliLeviWinternitz08a}:
\begin{equation}
  \label{eq:QQ-antiderivative}
  [\PP f](t) = -\half t \int\limits_0^t f(s)\,\rmd s  \qquad \text{ for } f\in \FSpace{E}{2}.
\end{equation}
It is known as a connection between \(\rme^{\QQ}\) and the (Gauss--)Weierstrass transform~\cite{Bilodeau62a}. Therefore, exponentiation of \(\QQ\) leads to the diffusion semigroup:
\begin{displaymath}
  \rme^{u\QQ} f(t) = \frac{1}{2 \sqrt{\pi u}}\int\limits_{-\infty}^\infty f(s)\, \rme^{-(t-s)^2/(4u)}\,\rmd s.
\end{displaymath}
Since the functional \(l_0\) is the evaluation \(\scalar[*]{l_0}{f}= f(0)\),  the covariant transform is:
\begin{equation}
  \label{eq:covariant-D2}
  \tilde{f}(u) =\scalar[*]{l_0}{\rme^{u \QQ} f} = \frac{1}{2\sqrt{\pi u}}\int\limits_{-\infty}^\infty f(s)\, \rme^{-s^2/(4u)}\,\rmd s.
\end{equation}
This formula also follows from the observation that the intertwining property~\eqref{eq:intertwining} between \(\QQ=\frac{\rmd^2\ }{\rmd t^2}\) and \(\frac{\rmd\ }{\rmd u}\) requires the integral kernel \(k(u,t)=\frac{1}{2\sqrt{\pi u}}\rme^{- t^2/(4u)} \), which is the fundamental solution of the heat equation \(\frac{\rmd k}{\rmd u}   = \frac{\rmd^2 k}{\rmd t^2}\).

We can use the same approach to find the intertwining operator with the momentum representation, cf. \S\ref{sec:fourier-transform}, to avoid a challenging exponentiation of the antiderivative operator \eqref{eq:QQ-antiderivative} as required by~\eqref{eq:intertwining-second}.  To intertwine operator \(\QQ=\frac{\rmd^2\ }{\rmd t^2}\) and \(vI\)  a kernel \(k(v,t)\) shall satisfy to the equation \(vk   = \frac{\rmd^2 k}{\rmd t^2}\) with the initial values \(k(v,0)=1\) and \(k_t'(v,0)=0\). The same answer appears as the generating function~\eqref{eq:generating-function-archetypal}:
\begin{equation}
  \label{eq:second-derivative-transform}
  k(v, t) =  \cos\left( \sqrt{v} t \right) = \sum_{k=0}^\infty \frac{(-v)^k\, t^{2k}}{(2k)!}
  \quad \text{ and the transform is } \quad
  \hat{f}(v) =\int\limits_{-\infty}^\infty f(t)  \cos\left( \sqrt{v} t \right) \rmd t\, .
\end{equation}
The operator \(\QQ =  -\frac{\rmd^2\ }{\rmd t^2} \) can be viewed as an ``umbra'' for a generic Schr\"odinger operator \(\frac{\rmd^2\ }{\rmd t^2} + V(t)I\), see~\cite{GelfandLevitan51a} for a related discussion with a usage of the orthogonalisation~\eqref{eq:orthogonal-polynomials}.


\subsubsection{Bessel operator and umbral calculus}
\label{sec:bessel-umbral}

Let \(\QQ=\oper{B}_\nu \) is the singular Bessel differential operator~\eqref{eq:bessel-operator}.
The umbral sequence is formed by even-order monomials with a tailored scaling, cf.~\eqref{eq:QQ-antiderivative}:
\begin{equation}
  \label{eq:bessel-umbral-sequence}
  p_n(t)=\frac{t^{2n}}{4^n \, n!\,\Gamma(n+\frac{\nu-1}{2})}\quad \text {with} \quad
  \Raise t^{2n} = \frac{ t^{2(n+1)}}{2(2n+\nu+1)} \quad \text{or} \quad
 [\Raise f](t)  = \half s^{1-\nu} \int\limits_0^t s^{\nu}\cdot f(s)\,\rmd s  .
\end{equation}
 The umbral functional is again the evaluation \(\scalar[*]{l_0}{f} =  f(0) \).
 Alternatively, an umbral sequence is formed by fractional powers:
\begin{displaymath}
  p_n(t)=\frac{t^{2n+1-\nu}}{4^n \, n!\,\Gamma(n+\frac{1-\nu}{2})}\quad \text {with} \quad
  \Raise t^{2n+1-\nu} = \frac{ t^{2n+3-\nu}}{2(2n+3-\nu)} \quad \text{or} \quad
 [\Raise f](t)  = \half  \int_0^t s\cdot f(s)\,\rmd s  .
 \end{displaymath}
 In this case the umbral functional is \(\scalar[*]{l_0}{f} = \lim_{t\rightarrow 0} (t^{\nu-1} f(t) )\).

 We are continuing with the first umbral sequence~\eqref{eq:bessel-umbral-sequence}.
 The generating function~\eqref{eq:generating-function-archetypal} is also a solution of the eigenvalue problem~\eqref{eq:L-eigenvalue-problem}, which is given by the normalised Bessel functions \(j_{\frac{\nu-1}{2}} (\sqrt{s}  t)\)~\eqref{eq:little-Bessel-function}:
 \begin{equation}
\label{eq:bessel-generating}
   j_{\frac{\nu-1}{2}} (\sqrt{s}  t) 
   =   \sum_{k=0}^\infty \frac{(-s)^k\, t^{2k}}{4^k \, k!\,\Gamma(k+\frac{1-\nu}{2})}.
 \end{equation}
 A pairing invariant under the action \(  \rme^{y\oper{B}_\nu}\) for all real \(y\) is
 \begin{displaymath}
   \scalar{g}{f} =  \int\limits_0^\infty g(t) \, f(t)\, t^\nu\,\rmd t.
 \end{displaymath}
 The respective induced covariant transform \(\scalar{j_{\frac{\nu-1}{2}} (\sqrt{s}  t)}{f(t)}\) is the Hankel transform~\eqref{eq:hankel-transform}.

 Thereafter, we can define a generalised translation \(T_t^y\)~\eqref{eq:generalised-translation} for the Bessel operator \(\oper{B}_\nu\).
 Then, the generating function \( j_{\frac{\nu-1}{2}} (\sqrt{s}  t)\) has the character property~\eqref{eq:generalised-character-property}, which takes the form~\eqref{eq:generalised-translation-defn}. Since the Hankel transform~\eqref{eq:hankel-transform} represents elements of \(\FSpace{E}{}\) as superpositions of the normalised Bessel functions we can write an integral representation of the generalised translations \(T_t^y\). Finally, the Poisson operator~\eqref{eq:poisson-operator} is the intertwining operator~\eqref{eq:poison-intertwining} between two umbral schemes generated by the Bessel operator and the second derivative considered in the previous paragraph.

\subsection{PDO type calculus of operators}
\label{sec:pdo-type-calculus}

An important tool to study boundary value problems with the Bessel operator \(\oper{B}_\nu\)~\eqref{eq:bessel-operator} is a \emph{generalised convolution}%
\index{generalised!convolution}%
\index{convolution!generalised} \cite[\S1.3]{Kipriyanov97a}, \cite[\S2.2.3]{SitnikShishkina19a}, \cite[\S3.4.3]{ShishkinaSitnik20a}. It can be defined through  the Hankel transform \(\oper{H}_\nu\)~\eqref{eq:hankel-transform} treated as a sort of  the Fourier transforms:
\begin{displaymath}
  \oper{H}_\nu(f*g) = \oper{H}_\nu(f) \cdot \oper{H}_\nu(g).
\end{displaymath}
Equivalently, the operator of convolution  \(\oper{C}_f: g \mapsto f*g\) is an integral operator of the generalised translations~\eqref{eq:generalised-translation-defn}:
\begin{displaymath}
  [\oper{C}_fg](t) = \int\limits_0^\infty f(y) \, [T_t^yg](t) \, y^\nu \,\rmd y  =
  \int\limits_0^\infty f(y) \, \sum_{k=0}^{\infty}  \phi_k(y)\,   [\oper{B}_\nu^k  g](t) \, y^\nu \,\rmd y.
\end{displaymath}
Clearly, this line of thoughts is applicable to any other umbral framework as well. From general principles, a large class of operators commuting with \(\Lower\) is contained in a weak closure of powers \(\Lower^k\), cf. Prop.~\ref{pr:binomial-type}.
Yet, we can use the umbral calculus to get a larger class of operators on \(\FSpace{E}{}\) as integrated representations of \(\uir{}{\FSpace{E}{}}\)~\eqref{eq:H-representation}. That is, for a suitable kernel \(k(x,y)\) we define the operator of \emph{relative convolution}%
\index{relative convolution}%
\index{convolution!relative }~\cite{Kisil94e,Kisil13a}:
\begin{equation}
  \label{eq:integrated-representation}
  \uir{}{}(k) = \int\!\!\int k(x,y)\, \rme^{y \QQ} \, \rme^{x \PP}\,\rmd x\,\rmd y.
\end{equation}
Of course, the expression is only meaningful if the operator \(\PP\) is defined at all, it adds another reason for consideration of \(\Raise\) along with \(\Lower\) to already presented in Rem.~\ref{re:role-Heisenberg-group}. The corresponding calculus of operators is well known~\cite{Howe80b,Folland89,Kisil19b}. In particular, the composition formula for two operators is
\begin{displaymath}
  \uir{}{}(k_1)\uir{}{}(k_2)= \uir{}{}(k_1 \natural k_2), \quad \text{where} \quad
  [k_1 \natural k_2](x,y)= \int\!\!\int k_1(x',y')\, k_2(x-x',y-y')\, \rme^{\pi \rmi(x'y-y'x)}\,\rmd x'\,\rmd y',
\end{displaymath}
and  \(k_1 \natural k_2\) is called the \emph{twisted convolution}%
\index{twisted convolution}%
\index{convolution!twisted} for the Heisenberg group. A further facilitation is provided through considering operators of the form \(a(\Raise,\Lower)\coloneqq\uir{}{}(\hat{a})\)~\eqref{eq:integrated-representation} with a kernel being the Fourier transform of a function \(a(x,y)\). If \(\uir{}{}\) is the Schr\"odinger representation~\eqref{eq:schrodinger-repres} then the operator \(\uir{}{}(\hat{a})\) is a \emph{pseudodifferential operator} (PDO)%
\index{PDO|see{pseudodifferential operator}}%
\index{pseudodifferential operator (PDO)}%
\index{operator!pseudodifferential} \(a(\Raise, \Lower)\) with the \emph{symbol}%
\index{symbol!pseudodifferential operator}%
\index{operator!pseudodifferential!symbol} \(a(x, y)\)~\cite{Howe80b,Folland89,Kisil19b}.

For other umbral frameworks we obtain different realisations of operators~\eqref{eq:integrated-representation}. However, all of them enjoy the same set of properties which are inherited from PDO through the Umbral principle. 
To make this approach successful for a specific class of operators one needs to make a suitable choice of an umbral model. On one hand selected operators \(\Raise\) and \(\Lower\) need to be accessible themselves, on the other hand \(\Raise\) and \(\Lower\) shall be sufficiently versatile to represent the desired class of operators as PDO \(a(\Raise, \Lower)\) with a manageable type of symbols \(a(x,y)\). 
A similar dilemma is elaborated for a different PDO-like calculus of operators in~\cite{DelgadoRuzhanskyTokmagambetov17a,CardonaKumarRuzhanskyTokmagambetov21a}.

\section{Discussion and closing remarks}
\label{sec:dicuss-clos-remarks}

In this paper we present umbral calculus foundations through representations of the Heisenberg group. Although this approach was around for a while~\cite{Cigler78,Feinsilver88a,Kisil97b,Kisil98d,Kisil01b,Kwasniewski03a,LeviTempestaWinternitz04a,DattoliLeviWinternitz08a,FaustinoRen13a} it largely remains out of the mainstream theory~\cite{RotaWay,StanleyI}. The discussed viewpoint on the umbral techniques covers some additional areas, e.g. generalised translations~\cite{Kipriyanov97a,Levitan49a,Levitan73a,SitnikShishkina19a,ShishkinaSitnik20a}. Such inclusion rewards the umbral approach with a removal of some unessential limitations, e.g. on the degree of a delta operator. Thus, it is stimulating to consider the umbral calculus not just as a one-way road to replicate some results from the archetypal model of power series to other situations. A more fruitful umbral ideology creates a hub, which facilities all-to-all exchanges between various umbral implementations, schematically depicted as follows:
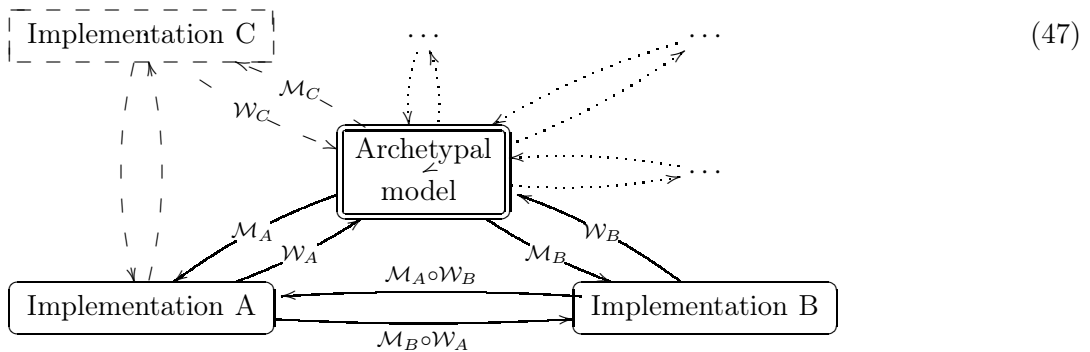
\begin{equation}
  \label{eq:umbral-hub}
  \xymatrix{ *+[F--:<3pt>]{\txt{ Implementation C \strut}}  \ar@{-->}@/_.65pc/[dr]|-{\oper{W}_C} \ar@/_.65pc/@{-->}[dd] & \cdots \ar@/_/@{.>}[d]& \cdots \ar@/_/@{.>}[dl]\\
    & *+[F=:<3pt>]{\txt{ Archetypal\ \strut\\model \strut}}  \ar\ar@{-->}@/_.65pc/[ul]|-{\oper{M}_C} \ar@/_.65pc/[dl]|-{\oper{M}_A} \ar@/_.65pc/[dr]|-{\oper{M}_B} \ar@/_/@{.>}[u] \ar@/_/@{.>}[ur] \ar@/_/@{.>}[r]& \cdots \ar@/_/@{.>}[l]\\
    *+[F-:<3pt>]{\txt{ Implementation A \strut}} \ar@/_.65pc/[ur]|-{\oper{W}_A} \ar@/_/[rr]_{\oper{M}_B\circ\oper{W}_A} \ar@/_.65pc/@{-->}[uu]& \qquad &*+[F-:<3pt>]{\txt{ Implementation B \strut}} \ar@/_.65pc./[ul]|-{\oper{W}_B} \ar@/_/[ll]_{\oper{M}_A\circ\oper{W}_B} &
    }
\end{equation}
Here the contravariant transform \(\oper{M}_A: \FSpace{\underline{E}}{} \rightarrow \FSpace{E}{A}\) is the adjoint of the covariant transform \(\oper{W}_A: \FSpace{E}{A} \rightarrow \FSpace{\underline{E}}{}\). Then, the transmutation \(\oper{T}_{BA} = \oper{M}_B\circ\oper{W}_A: \FSpace{E}{A} \rightarrow \FSpace{E}{B}\) intertwines respective ladder operator
:
\begin{displaymath}
  \oper{T}_{BA} \circ \Lower_A= \Lower_B  \circ \oper{T}_{BA}  \qquad \text{ and } \qquad
  \oper{T}_{BA} \circ \Raise_A= \Raise_B  \circ \oper{T}_{BA} .
\end{displaymath}
Such construction of transmutations falls within the \emph{composition method}~\cite[Ch.~5]{SitnikShishkina19a}, \cite[Ch.~6]{ShishkinaSitnik20a}.

Now we can formulate the extension of the original Umbral Principle.

\begin{principle}[Extended]
  For any two implementations of the umbral model intertwining operators (transmutations) \(\oper{M}_B\circ\oper{W}_A\) and \(\oper{M}_A\circ\oper{W}_B\)~\eqref{eq:umbral-hub} allow us to exchange  scopes, problems, ideas, methods, results, etc. between those implementations.     
\end{principle}
Here are few examples of such cross-fertilisation:
\begin{itemize}
\item The generalised translations for the Bessel operators shows that the umbral delta operator shall not be restricted by the condition that it reduces the order of polynomials by one. It will be interesting to see if the Bessel operator~\eqref{eq:bessel-operator} may be useful as a delta in some combinatorial and enumeration problems.
\item The generalised translations \(T^y\) naturally create respective convolutions. Yet, we can transport a technique for calculus of larger classes (not necessarily \(T^y\)-invariant) operators from the theory of PDO~\cite{Howe80b,Folland89,Kisil19b} generated by the archetypal umbral model from Ex.~\ref{ex:archetypal-polinom}.
\item It is possible to transport the idea of the Feynman path integral from quantum mechanics to combinatorics using the umbral framework~\cite{Kisil98d}. It may be interesting to see the meaning  of path integrals for generalised translations as well as other umbral implementations.
\end{itemize}

Despite of the exceptional r\^ole of the Heisenberg group in the analysis~\cite{Howe80a} an umbral approach may be build on other group representations as well. The next natural candidate to implement a sort of ladder actions~\eqref{eq:hat-Q}--\eqref{eq:hat-P} would be the group \(\SL\). Indeed, for a given basis \(p_n\) of a linear space \(\FSpace{E}{}\) there are many essentially different possibilities to chose sequences of numbers \({a_n}\), \({b_n}\), \({c_n}\) such that the extended set \(\{\Lower, \Raise, \oper{Z}\}\) of linear operators defined by~\cite{HoweTan92,Kisil11a}
\begin{equation}
  \label{eq:sl2ladder-action}
  \Lower p_n = a_n p_{n-1}, \qquad \Raise p_n = b_n p_{n+1}, \qquad \oper{Z} p_n = c_n p_n \qquad
\end{equation}
make a representation of the Lie algebra of \(\SL\):
\begin{equation}
  \label{eq:sl2-lie-algebra}
  [\Lower,\Raise]= \lambda \oper{Z} , \qquad[\oper{Z},\Lower]= \lambda_- \Lower, \qquad  [\oper{Z},\Raise]= \lambda_+ \Raise, \qquad \text{for some scalars} \quad \lambda, \ \lambda_-, \ \lambda_+.
\end{equation}
A graphical illustration of action~\eqref{eq:sl2ladder-action} is:
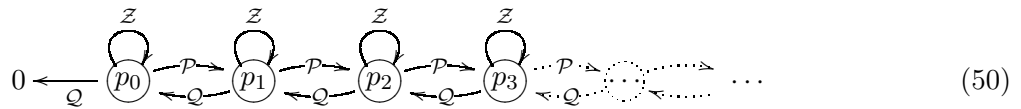
\begin{equation}
  \label{eq:sl2-action-illutration}
  \xymatrix{ 0 
    &*+[o][F]{\ p_0\ } \ar@/^/[r]|-{\Raise} \ar[l]^{\Lower}     \ar@(ul,ur)^{\oper{Z}}[] 
    &*+[o][F]{\ p_1\ }\ar@/^/[r]|-{\Raise} \ar@/^/[l]|-{\Lower}
    \ar@(ul,ur)^{\oper{Z}}[] 
    &*+[o][F]{\ p_2\ }\ar@/^/[r]|-{\Raise} \ar@/^/[l]|-{\Lower}
    \ar@(ul,ur)^{\oper{Z}}[]
    &*+[o][F]{\ p_3\ }\ar@/^/@{..>}[r]|-{\Raise} \ar@/^/[l]|-{\Lower}
    \ar@(ul,ur)^{\oper{Z}}[]
    &*+[o][F.]{\cdots}\ar@/^/@{..>}[l]|-{\Lower}     \ar@/^/@{..>}[r]
     &\ \cdots\ 
     \ar@/^/@{..>}[l]
   } 
\end{equation}
where, in comparison to~\eqref{eq:ladder-action},  the reflexive arrows for \(\oper{Z}\)-action are added. 
Such ladders include much more than just the unitary representations of \(\SL\) by the discrete holomorphic series~\cite{HoweTan92,Kisil11a}. The additional representations have various applications, e.g. for quasi-exact solvable quantum systems~\cite{Turbiner88a}. 

Interestingly, the Heisenberg ladder~\eqref{eq:ladder-action} can be a source of some \(\SL\) representations of the form~\eqref{eq:sl2-action-illutration} through the  quadratic algebra concept~\cite[\textsection{}2.2.4]{Gazeau09a}. Indeed, if \([\Raise,\Lower]=-\iota I\) then operators \(\Raise_2= \Raise^2\), \(\Lower_2=\Lower^2\) and \(\oper{Z}_2= \Raise\Lower\) satisfies commutators~\eqref{eq:sl2-lie-algebra} and form a \emph{metaplectic representation}%
\index{metaplectic representation}%
\index{representation!metaplectic} of \(\SL\)~\cite{Folland89,HoweTan92}. Therefore, the even-numbered elements of the ladder~\eqref{eq:ladder-action} becomes nodes of the action~\eqref{eq:sl2-action-illutration} depicted as follows:
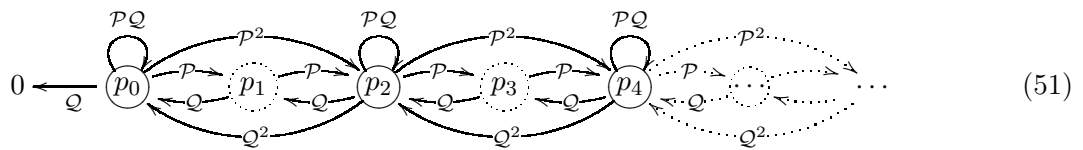
\begin{equation}
  \label{eq:factorised-action}
  \xymatrix{ 0 
    &*+[o][F]{\ p_0\ } \ar@/^/[r]|-{\Raise} \ar@/_/@(ur,ul)[rr]|-{\Raise^2} \ar[l]^{\Lower}     \ar@(ul,ur)^{\Raise\Lower}[] 
    &*+[o][F.]{\ p_1\  }\ar@/^/[r]|-{\Raise} \ar@/^/[l]|-{\Lower} 
    &*+[o][F]{\ p_2\ }\ar@/^/[r]|-{\Raise} \ar@/^/[l]|-{\Lower}
    \ar@/_/@(ur,ul)[rr]|-{\Raise^2} \ar@/_/@(dl,dr)[ll]|-{\Lower^2}
    \ar@(ul,ur)^{\Raise\Lower}[] 
    &*+[o][F.]{\ p_3\ } \ar@/^/[r]|-{\Raise} \ar@/^/[l]|-{\Lower} 
    &*+[o][F]{\ p_4\ }\ar@/^/@{..>}[r]|-{\Raise} \ar@/^/[l]|-{\Lower}
        \ar@/_/@{..>}@(ur,ul)[rr]|-{\Raise^2} \ar@/_/@(dl,dr)[ll]|-{\Lower^2}     \ar@(ul,ur)^{\Raise\Lower}[]
    &*+[o][F.]{\cdots}\ar@/^/@{..>}[l]|-{\Lower}     \ar@/^/@{..>}[r]
     &\ \cdots\ 
    \ar@/^/@{..>}[l]
    \ar@/_/@{..>}@(dl,dr)[ll]|-{\Lower^2} 
  }
\end{equation}
A transition in the opposite direction---from the \(\SL\)-action~\eqref{eq:sl2-action-illutration} to its embedding in an extended ladder~\eqref{eq:factorised-action}---can be viewed as a sort of factorisation of a second-order operator. This provides a more detailed resolutions, say, for the Bessel operator \(\oper{B}_\nu\)~\eqref{eq:bessel-operator}.

Overall, the umbral approach in the context of \(\SL\) may be an exciting topic with an extensive exchange of ideas between various fields and deserves a further study. 



\let\germ=\mathfrak
\newfam\cyrfam
\font\tencyr=wncyr10 \font\sevencyr=wncyr7 \font\fivecyr=wncyr5
\def\cyr{\fam\cyrfam\tencyr\cyracc}
\input{cyracc.def}
\textfont\cyrfam=\tencyr \scriptfont\cyrfam=\sevencyr \scriptscriptfont\cyrfam=\fivecyr


\begin{thebibliography}{47}
\providecommand{\natexlab}[1]{#1}
\providecommand{\url}[1]{\texttt{#1}}
\expandafter\ifx\csname urlstyle\endcsname\relax
  \providecommand{\doi}[1]{doi: #1}\else
  \providecommand{\doi}{doi: \begingroup \urlstyle{rm}\Url}\fi

\bibitem[Kipriyanov(1997)]{Kipriyanov97a}
I.~A. Kipriyanov.
\newblock \emph{Singular Boundary Value Problems \cyr [{S}ingulyarnye
  \`Ellipticheskie Kraevye Zadachi]}.
\newblock Fizmatlit ``Nauka'', Moscow, 1997.
\newblock ISBN 5-02-014799-0.
\newblock In Russian.

\bibitem[Sitnik and Shishkina(2019)]{SitnikShishkina19a}
S.~Sitnik and E.~Shishkina.
\newblock \emph{Method of Transmutations for Differential Equations with
  {Bessel} Operators}.
\newblock Fizmatlit, Moscow, 2019.
\newblock (in Russian).

\bibitem[Shishkina and Sitnik(2020)]{ShishkinaSitnik20a}
E.~Shishkina and S.~Sitnik.
\newblock \emph{Transmutations, Singular and Fractional Differential Equations
  with Applications to Mathematical Physics}.
\newblock Mathematics in Science and Engineering. Elsevier/Academic Press,
  London, 2020.
\newblock ISBN 978-0-12-819781-3.

\bibitem[Kravchenko and Sitnik(2020)]{KravchenkoSitnik20a}
V.~V. Kravchenko and S.~M. Sitnik.
\newblock Some recent developments in the transmutation operator approach.
\newblock In V.~V. Kravchenko and S.~M. Sitnik, editors, \emph{Transmutation
  Operators and Applications}, pages 3--9. Springer International Publishing,
  Cham, 2020.
\newblock ISBN 978-3-030-35914-0.
\newblock URL \url{https://doi.org/10.1007/978-3-030-35914-0_1}.

\bibitem[Levitan(1949)]{Levitan49a}
B.~M. Levitan.
\newblock The application of generalized displacement operators to linear
  differential equations of the second order.
\newblock \emph{Uspehi Matem. Nauk (N.S.)}, 4\penalty0 (1(29)):\penalty0
  3--112, 1949.
\newblock ISSN 0042-1316.

\bibitem[Levitan(1951)]{Levitan51a}
B.~M. Levitan.
\newblock Expansion in {F}ourier series and integrals with {B}essel functions.
\newblock \emph{Uspehi Matem. Nauk (N.S.)}, 6\penalty0 (2(42)):\penalty0
  102--143, 1951.
\newblock ISSN 0042-1316.

\bibitem[Levitan(1973)]{Levitan73a}
B.~M. Levitan.
\newblock \emph{{\cyr Teoriya operatorov obobshchennogo sdviga}}.
\newblock Izdat. ``Nauka'', Moscow, 1973.

\bibitem[Rota(1975)]{Rota75}
G.-C. Rota.
\newblock \emph{Finite Operator Calculus}.
\newblock Academic Press, Inc., New York, 1975.

\bibitem[Kung et~al.(2009)Kung, Rota, and Yan]{RotaWay}
J.~P.~S. Kung, G.-C. Rota, and C.~H. Yan.
\newblock \emph{Combinatorics: the {R}ota way}.
\newblock Cambridge Mathematical Library. Cambridge University Press,
  Cambridge, 2009.
\newblock ISBN 978-0-521-73794-4.

\bibitem[Howe(1980{\natexlab{a}})]{Howe80a}
R.~Howe.
\newblock On the role of the {Heisenberg} group in harmonic analysis.
\newblock \emph{Bull. Amer. Math. Soc. (N.S.)}, 3\penalty0 (2):\penalty0
  821--843, 1980{\natexlab{a}}.
\newblock ISSN 0002-9904.

\bibitem[Kisil(1996)]{Kisil95i}
V.~V. Kisil.
\newblock M\"obius transformations and monogenic functional calculus.
\newblock \emph{Electron. Res. Announc. Amer. Math. Soc.}, 2\penalty0
  (1):\penalty0 26--33, 1996.
\newblock ISSN 1079-6762.
\newblock
  \href{http://www.ams.org/era/1996-02-01/S1079-6762-96-00004-2/}{On-line}.

\bibitem[Kisil(2004)]{Kisil02a}
V.~V. Kisil.
\newblock Spectrum as the support of functional calculus.
\newblock In \emph{Functional analysis and its applications}, volume 197 of
  \emph{North-Holland Math. Stud.}, pages 133--141, Amsterdam, 2004. Elsevier.
\newblock \arXiv{math.FA/0208249}.

\bibitem[Kisil(2012{\natexlab{a}})]{Kisil12b}
V.~V. Kisil.
\newblock Operator covariant transform and local principle.
\newblock \emph{J. Phys. A: Math. Theor.}, 45:\penalty0 244022,
  2012{\natexlab{a}}.
\newblock \doi{10.1088/1751-8113/45/24/244022}.
\newblock \arXiv{1201.1749}.
  \href{http://stacks.iop.org/1751-8121/45/244022}{On-line}.

\bibitem[Ali et~al.(2000)Ali, Antoine, and Gazeau]{AliAntGaz00}
S.~T. Ali, J.-P. Antoine, and J.-P. Gazeau.
\newblock \emph{Coherent States, Wavelets and Their Generalizations}.
\newblock Graduate Texts in Contemporary Physics. Springer-Verlag, New York,
  2000.
\newblock ISBN 0-387-98908-0.

\bibitem[Perelomov(1986)]{Perelomov86}
A.~Perelomov.
\newblock \emph{Generalized Coherent States and Their Applications}.
\newblock Texts and Monographs in Physics. Springer-Verlag, Berlin, 1986.
\newblock ISBN 3-540-15912-6.

\bibitem[Kisil(2012{\natexlab{b}})]{Kisil11c}
V.~V. Kisil.
\newblock {E}rlangen programme at large: {An} overview.
\newblock In S.~Rogosin and A.~Koroleva, editors, \emph{Advances in Applied
  Analysis}, chapter~1, pages 1--94. Birkh\"auser Verlag, Basel,
  2012{\natexlab{b}}.
\newblock \arXiv{1106.1686}.

\bibitem[Feichtinger and Gr\"{o}chenig(1988)]{FeichtingerGrochenig88a}
H.~G. Feichtinger and K.~Gr\"{o}chenig.
\newblock A unified approach to atomic decompositions via integrable group
  representations.
\newblock In \emph{Function spaces and applications ({L}und, 1986)}, volume
  1302 of \emph{Lecture Notes in Math.}, pages 52--73. Springer, Berlin, 1988.
\newblock \doi{10.1007/BFb0078863}.
\newblock URL \url{https://doi.org/10.1007/BFb0078863}.

\bibitem[Kisil(1999{\natexlab{a}})]{Kisil98a}
V.~V. Kisil.
\newblock Wavelets in {B}anach spaces.
\newblock \emph{Acta Appl. Math.}, 59\penalty0 (1):\penalty0 79--109,
  1999{\natexlab{a}}.
\newblock ISSN 0167-8019.
\newblock \arXiv{math/9807141},
  \href{http://dx.doi.org/10.1023/A:1006394832290}{On-line}.

\bibitem[Folland(1989)]{Folland89}
G.~B. Folland.
\newblock \emph{Harmonic Analysis in Phase Space}, volume 122 of \emph{Annals
  of Mathematics Studies}.
\newblock Princeton University Press, Princeton, NJ, 1989.
\newblock ISBN 0-691-08527-7; 0-691-08528-5.

\bibitem[Feinsilver(1988)]{Feinsilver88a}
P.~Feinsilver.
\newblock Lie algebras and recurrence relations. {I}.
\newblock \emph{Acta Appl. Math.}, 13\penalty0 (3):\penalty0 291--333, 1988.
\newblock ISSN 0167-8019.
\newblock \doi{10.1007/BF00046967}.
\newblock URL \url{https://doi.org/10.1007/BF00046967}.

\bibitem[Kisil(2012{\natexlab{c}})]{Kisil10a}
V.~V. Kisil.
\newblock Hypercomplex representations of the {H}eisenberg group and mechanics.
\newblock \emph{Internat. J. Theoret. Phys.}, 51\penalty0 (3):\penalty0
  964--984, 2012{\natexlab{c}}.
\newblock ISSN 0020-7748.
\newblock \doi{10.1007/s10773-011-0970-0}.
\newblock URL \url{http://dx.doi.org/10.1007/s10773-011-0970-0}.
\newblock \arXiv{1005.5057}. \Zbl{1247.81232}.

\bibitem[Al~Ameer and Kisil(2022)]{AlameerKisil21a}
A.~A. Al~Ameer and V.~V. Kisil.
\newblock Tuning co- and contra-variant transforms: the {Heisenberg} group
  illustration.
\newblock \emph{SIGMA}, 18\penalty0 (065):\penalty0 21 p., 2022.
\newblock \arXiv{2105.13811}.

\bibitem[Kisil(2014{\natexlab{a}})]{Kisil12d}
V.~V. Kisil.
\newblock The real and complex techniques in harmonic analysis from the point
  of view of covariant transform.
\newblock \emph{Eurasian Math. J.}, 5:\penalty0 95--121, 2014{\natexlab{a}}.
\newblock \arXiv{1209.5072}.
  \href{http://emj.enu.kz/images/pdf/2014/5-1-4.pdf}{On-line}.

\bibitem[Kisil(2010)]{Kisil09d}
V.~V. Kisil.
\newblock Wavelets beyond admissibility.
\newblock In M.~Ruzhansky and J.~Wirth, editors, \emph{Progress in Analysis and
  its Applications}, pages 219--225. World Sci. Publ., Hackensack, NJ, 2010.
\newblock URL \url{http://dx.doi.org/10.1142/9789814313179_0029}.
\newblock \arXiv{0911.4701}. \Zbl{1269.30052}.

\bibitem[Kisil(2014{\natexlab{b}})]{Kisil13a}
V.~V. Kisil.
\newblock Calculus of operators: {C}ovariant transform and relative
  convolutions.
\newblock \emph{Banach J. Math. Anal.}, 8\penalty0 (2):\penalty0 156--184,
  2014{\natexlab{b}}.
\newblock URL \url{http://www.emis.de/journals/BJMA/tex_v8_n2_a15.pdf}.
\newblock \arXiv{1304.2792},
  \href{http://www.emis.de/journals/BJMA/tex_v8_n2_a15.pdf}{on-line}.

\bibitem[Kung(1995)]{Rota95}
J.~P. Kung, editor.
\newblock \emph{Gian-Carlo Rota on Combinatorics: Introductory Papers and
  Commentaries}, volume~1 of \emph{Contemporary Mathematicians}.
\newblock Birkh\"auser Verlag, Boston, 1995.

\bibitem[Cigler(1978)]{Cigler78}
J.~Cigler.
\newblock Some remarks on {Rota's} umbral calculus.
\newblock \emph{Nederl. Akad. Wetensch. Proc. Ser. A}, 81:\penalty0 27--42,
  1978.

\bibitem[Kisil(2000)]{Kisil97b}
V.~V. Kisil.
\newblock Umbral calculus and cancellative semigroup algebras.
\newblock \emph{Z. Anal. Anwendungen}, 19\penalty0 (2):\penalty0 315--338,
  2000.
\newblock ISSN 0232-2064.
\newblock \arXiv{funct-an/9704001}. \Zbl{0959.43004}.

\bibitem[Kisil(2002{\natexlab{a}})]{Kisil01b}
V.~V. Kisil.
\newblock Tokens: an algebraic construction common in combinatorics, analysis,
  and physics.
\newblock In \emph{Ukrainian Mathematics Congress---2001 (Ukrainian)}, pages
  146--155. Nats\=\i onal. Akad. Nauk Ukra\"\i ni \=Inst. Mat., Kiev,
  2002{\natexlab{a}}.
\newblock \arXiv{math.FA/0201012}.

\bibitem[Kwa\'{s}niewski(2003)]{Kwasniewski03a}
A.~K. Kwa\'{s}niewski.
\newblock Main theorems of extended finite operator calculus.
\newblock \emph{Integral Transforms Spec. Funct.}, 14\penalty0 (6):\penalty0
  499--516, 2003.
\newblock ISSN 1065-2469.
\newblock \doi{10.1080/10652460290029743}.
\newblock URL \url{https://doi.org/10.1080/10652460290029743}.

\bibitem[Levi et~al.(2004)Levi, Tempesta, and
  Winternitz]{LeviTempestaWinternitz04a}
D.~Levi, P.~Tempesta, and P.~Winternitz.
\newblock Umbral calculus, difference equations and the discrete
  {S}chr\"{o}dinger equation.
\newblock \emph{J. Math. Phys.}, 45\penalty0 (11):\penalty0 4077--4105, 2004.
\newblock ISSN 0022-2488.
\newblock \doi{10.1063/1.1780612}.
\newblock URL \url{https://doi.org/10.1063/1.1780612}.

\bibitem[Dattoli et~al.(2008)Dattoli, Levi, and
  Winternitz]{DattoliLeviWinternitz08a}
G.~Dattoli, D.~Levi, and P.~Winternitz.
\newblock {H}eisenberg algebra, umbral calculus and orthogonal polynomials.
\newblock \emph{J. Math. Phys.}, 49\penalty0 (5):\penalty0 053509, 19, 2008.
\newblock ISSN 0022-2488.
\newblock \doi{10.1063/1.2909731}.
\newblock URL \url{https://doi.org/10.1063/1.2909731}.

\bibitem[Faustino and Ren(2011)]{FaustinoRen13a}
N.~Faustino and G.~Ren.
\newblock ({D}iscrete) {A}lmansi type decompositions: An umbral calculus
  framework based on {$\germ{osp}(1|2)$} symmetries.
\newblock \emph{Math. Methods Appl. Sci.}, 34\penalty0 (16):\penalty0
  1961--1979, 2011.
\newblock ISSN 0170-4214.
\newblock \doi{10.1002/mma.1498}.
\newblock URL \url{https://doi.org/10.1002/mma.1498}.

\bibitem[Mullin and Rota(1970)]{FoundationIII}
R.~Mullin and G.-C. Rota.
\newblock On the foundation of combinatorial theory ({III}): Theory of binomial
  enumeration.
\newblock In B.Harris, editor, \emph{Graph Theory and Its Applications}, pages
  167--213. Academic Press, Inc., New York, 1970.
\newblock Reprinted in~\cite[pp.~118--147]{Rota95}.

\bibitem[Bilodeau(1962)]{Bilodeau62a}
G.~G. Bilodeau.
\newblock The {W}eierstrass transform and {H}ermite polynomials.
\newblock \emph{Duke Math. J.}, 29:\penalty0 293--308, 1962.
\newblock ISSN 0012-7094.
\newblock URL \url{http://projecteuclid.org/euclid.dmj/1077470134}.

\bibitem[Gel\cprime~fand and Levitan(1951)]{GelfandLevitan51a}
I.~M. Gel\cprime~fand and B.~M. Levitan.
\newblock On the determination of a differential equation from its spectral
  function.
\newblock \emph{Izvestiya Akad. Nauk SSSR. Ser. Mat.}, 15:\penalty0 309--360,
  1951.
\newblock ISSN 0373-2436.

\bibitem[Kisil(1999{\natexlab{b}})]{Kisil94e}
V.~V. Kisil.
\newblock Relative convolutions. {I}. {P}roperties and applications.
\newblock \emph{Adv. Math.}, 147\penalty0 (1):\penalty0 35--73,
  1999{\natexlab{b}}.
\newblock ISSN 0001-8708.
\newblock \arXiv{funct-an/9410001},
  \href{http://www.idealibrary.com/links/doi/10.1006/aima.1999.1833}{On-line}.
  \Zbl{933.43004}.

\bibitem[Howe(1980{\natexlab{b}})]{Howe80b}
R.~Howe.
\newblock Quantum mechanics and partial differential equations.
\newblock \emph{J. Funct. Anal.}, 38\penalty0 (2):\penalty0 188--254,
  1980{\natexlab{b}}.
\newblock ISSN 0022-1236.

\bibitem[Kisil(2023)]{Kisil19b}
V.~V. Kisil.
\newblock Cross-{Toeplitz} operators on the {Fock}--{Segal}--{Bargmann} space
  and two-sided convolutions on the {Heisenberg} group.
\newblock \emph{Ann. Funct. Anal.}, 14\penalty0 (38), 2023.
\newblock to appear. \arXiv{2108.13710}. \doi{10.1007/s43034-022-00249-7}.

\bibitem[Delgado et~al.(2017)Delgado, Ruzhansky, and
  Tokmagambetov]{DelgadoRuzhanskyTokmagambetov17a}
J.~Delgado, M.~Ruzhansky, and N.~Tokmagambetov.
\newblock Schatten classes, nuclearity and nonharmonic analysis on compact
  manifolds with boundary.
\newblock \emph{J. Math. Pures Appl. (9)}, 107\penalty0 (6):\penalty0 758--783,
  2017.
\newblock ISSN 0021-7824.
\newblock \doi{10.1016/j.matpur.2016.10.005}.
\newblock URL \url{https://doi.org/10.1016/j.matpur.2016.10.005}.

\bibitem[Cardona et~al.(2021)Cardona, Kumar, Ruzhansky, and
  Tokmagambetov]{CardonaKumarRuzhanskyTokmagambetov21a}
D.~Cardona, V.~Kumar, M.~Ruzhansky, and N.~Tokmagambetov.
\newblock Global functional calculus, lower/upper bounds and evolution
  equations on manifolds with boundary, 2021.

\bibitem[Kisil(2002{\natexlab{b}})]{Kisil98d}
V.~V. Kisil.
\newblock Polynomial sequences of binomial type and path integrals.
\newblock \emph{Ann. Comb.}, 6\penalty0 (1):\penalty0 45--56,
  2002{\natexlab{b}}.
\newblock ISSN 0218-0006.
\newblock \arXiv{math/9808040} \Zbl{1009.05013}.

\bibitem[Stanley(1997)]{StanleyI}
R.~P. Stanley.
\newblock \emph{Enumerative Combinatorics. {V}ol. 1}.
\newblock Cambridge University Press, Cambridge, 1997.
\newblock ISBN 0-521-55309-1.
\newblock With a foreword by Gian-Carlo Rota, Corrected reprint of the 1986
  original.

\bibitem[Howe and Tan(1992)]{HoweTan92}
R.~Howe and E.-C. Tan.
\newblock \emph{Nonabelian Harmonic Analysis. {Applications of
  ${{\rm{S}}L}(2,{{\bf{R}}})$}}.
\newblock Springer-Verlag, New York, 1992.
\newblock ISBN 0-387-97768-6.

\bibitem[Kisil(2011)]{Kisil11a}
V.~V. Kisil.
\newblock {E}rlangen {P}rogramme at {L}arge 3.2: {L}adder operators in
  hypercomplex mechanics.
\newblock \emph{Acta Polytechnica}, 51\penalty0 (4):\penalty0 44--53, 2011.
\newblock \href{http://ctn.cvut.cz/ap/download.php?id=614}{on-line}.
  \arXiv{1103.1120}.

\bibitem[Turbiner(1988)]{Turbiner88a}
A.~V. Turbiner.
\newblock Quasi-exactly-solvable problems and {${\rm sl}(2)$} algebra.
\newblock \emph{Comm. Math. Phys.}, 118\penalty0 (3):\penalty0 467--474, 1988.
\newblock ISSN 0010-3616.
\newblock URL \url{http://projecteuclid.org/euclid.cmp/1104162094}.

\bibitem[Gazeau(2009)]{Gazeau09a}
J.-P. Gazeau.
\newblock \emph{{Coherent States in Quantum Physics}}.
\newblock Wiley-VCH Verlag, 2009.
\newblock ISBN 9783527407095.
\newblock \doi{10.1002/9783527628285}.

\end{thebibliography}

\providecommand{\noopsort}[1]{} \providecommand{\printfirst}[2]{#1}
  \providecommand{\singleletter}[1]{#1} \providecommand{\switchargs}[2]{#2#1}
  \providecommand{\irm}{\textup{I}} \providecommand{\iirm}{\textup{II}}
  \providecommand{\vrm}{\textup{V}} \providecommand{\cprime}{'}
  \providecommand{\eprint}[2]{\texttt{#2}}
  \providecommand{\myeprint}[2]{\texttt{#2}}
  \providecommand{\arXiv}[1]{\myeprint{http://arXiv.org/abs/#1}{arXiv:#1}}
  \providecommand{\doi}[1]{\href{http://dx.doi.org/#1}{doi:
  #1}}\providecommand{\CPP}{\texttt{C++}}
  \providecommand{\NoWEB}{\texttt{noweb}}
  \providecommand{\MetaPost}{\texttt{Meta}\-\texttt{Post}}
  \providecommand{\GiNaC}{\textsf{GiNaC}}
  \providecommand{\pyGiNaC}{\textsf{pyGiNaC}}
  \providecommand{\Asymptote}{\texttt{Asymptote}}

\end{document}